\documentclass[12pt]{article}
\usepackage{amsmath,euscript}
\usepackage{amssymb}
\usepackage{amsthm}
\usepackage{enumerate}
\usepackage{amsfonts}
\usepackage{comment}
\usepackage{colortbl}
\usepackage{a4wide}
\usepackage{amssymb,amsmath,array}

\newtheorem{thm}{Theorem}[section]
\newtheorem{lem}[thm]{Lemma}

\newtheorem{prop}[thm]{Proposition}
\newtheorem{cor}[thm]{Corollary}
\newtheorem{conj}[thm]{Conjecture}

\theoremstyle{defn}
\newtheorem*{defn}{Definition}

\newtheorem*{remark}{Remark}

\newtheorem*{notn}{Notation}
\newtheorem{quest}{Question}

\newcommand{\Pf} {\rm Pfaffian}
\newcommand{\Ca} {\mathcal{C}_{m,n}}
\newcommand{\PMC}{\mathcal{PM}(C)}
\newcommand{\PMCi}{\mathcal{PM}(C_i)}
\newcommand{\VC}{{\rm Vert}(C)}
\newcommand{\Vpi}{{\rm Vert}(\pi)}

\newcommand{\mC}{\mathcal{C}}

\newcommand{\gray}{\cellcolor[gray]{0.1} }

\newcommand{\gc}{ [ \hspace{-0.65mm} [}
\newcommand{\dc}{]  \hspace{-0.65mm} ]}

\newcommand{\ia}{i,\alpha}

\newcommand{\fract}{\mathrm{Frac}}

\newcommand{\hc}{\mathcal{H}}

\newcommand{\gkdim}{{\rm GKdim}}
\newcommand{\Tdeg}{{\rm Tdeg}}

\newcommand{\spec}{{\rm Spec}}

\newcommand{\prim}{{\rm Prim}}

\def\C{\mathbb{K}}

\def\co{{\mathcal O}}

\def\oqmn{\co_q(M_n)}

\def\oqmm13{\co_q(M_{1,3})}
\def\oqm2n{\co_q(M_{2,n})}
\def\oqmmn{\co_q(M_{m,n})}
\def\oqmnm{\co_q(M_{n,m})}
\def\oqsln{\co_q(SL_n)}

\def\ia{i,\alpha}

\input{xy}
\xyoption{all}

\begin{document}

\title{Dimension and enumeration of primitive ideals in quantum algebras}
\author{J. Bell\thanks{The first author thanks NSERC for its generous support}, S. Launois\thanks{This research was supported by a Marie Curie Intra-European
  Fellowship within the $6^{\mbox{th}}$ European Community Framework
  Programme held at the University of Edinburgh and by Leverhulme Research Interchange Grant F/00158/X}\;  and N. Nguyen\;
}
\date{}

\maketitle


\begin{abstract}
In this paper, we study the primitive ideals of quantum algebras supporting a rational torus action.  We first prove a quantum analogue of a theorem of Dixmier; namely, we show that the Gelfand-Kirillov dimension of primitive factors of various quantum algebras is always even.  Next we give a combinatorial criterion for a prime ideal that is invariant under the torus action to be primitive.  We use this criterion to obtain a formula for the number of primitive ideals in the algebra of $2\times n$ quantum matrices that are invariant under the action of the torus.
Roughly speaking, this can be thought of as giving an enumeration of the points that are invariant under the induced action of the torus in the ``variety of $2\times n$ quantum matrices''. \end{abstract}


\section*{Introduction}

This paper is concerned with the primitive ideals of certain quantum algebras, and in particular with the primitive ideals of the algebra $\oqmmn$ of generic quantum matrices. Since S.P. Smith famous lectures on ring theoretic aspects of quantum groups in 1989 (see \cite{paul}), 
primitive ideals of quantum algebras have been extensively studied (see for instance \cite{bgbook} and \cite{josbook}). In particular, Hodges and Levasseur \cite{hl1,hl2} have discovered a remarkable partition of the primitive spectrum 
of the quantum special linear group $\oqsln$ and proved that the primitive ideals of $\oqsln$ correspond bijectively 
to the symplectic leaves in $SL_n$ (endowed with the semi-classical Poisson structure coming from the commutators of $\oqsln$). These results were next extended by Joseph to the standard quantised coordinate ring $\co_q(G)$ of a complex semisimple algebraic group $G$ \cite{jos1,jos2}. 
Let us mention however that it is not known (except in the case where $G=SL_2$) whether such bijection can be made into an homeomorphism. 
Later, it was observed by Brown and Goodearl \cite{bgtrans} that the existence of such partition relies very much on the action of a torus of automorphisms, and a general theory was then developed by Goodearl and Letzter in order to study the primitive spectrum of an algebra supporting a "nice" torus action \cite{gl2}. In particular, they constructed a partition, called the $H$-stratification, of the prime spectrum of such algebras which also induces by restriction a partition of the primitive spectrum of such algebras. This theory can be applied to many quantum algebras 
in the generic case, and in particular to the algebra $\oqmmn$ of generic quantum matrices as there is a natural action of the algebraic torus $\hc:= \mathbb{K}^{m+n}$ on this algebra.  In this case, the $H$-stratification theory of Goodearl and Letzter predicts the following \cite{gl2} (see also \cite{bgbook}). First, the number of prime ideals of $\oqmmn$ invariant under the action of this torus $\hc$ is finite. Next, the prime spectrum of $\oqmmn$ admits a stratification into finitely many $\hc$-strata. Each $\hc$-stratum is defined by a unique $\hc$-invariant prime ideal---that is minimal in its $\hc$-stratum---and is homeomorphic to the scheme of irreducible subvarieties of a torus. Moreover the primitive ideals correspond to those primes that are maximal in their $\hc$-strata and the Dixmier-Moeglin Equivalence holds. 

The first aim of this paper is to develop a strategy to recognise those $\hc$-invariant prime ideals that are primitive. In particular, we give a combinatorial criterion for an $\hc$-invariant prime ideal to be primitive. This generalises a result of Lenagan and the second author \cite{laulen} who gave a criterion for $(0)$ to be primitive. Our criterion in this paper is expressed in terms of combinatorial tools such as Cauchon diagrams---recently, Cauchon diagrams have also appeared in the literature under the name ''Le-diagrams'', see for instance \cite{p,w}---, perfect matchings, and Pfaffians of $0,\pm 1$ matrices.  We discuss these concepts in Sections 2.2 and 2.3.  As a corollary, we obtain a formula for the total number of  primitive $\hc$-invariant ideals in $\oqm2n$. 
More precisely, we show that the number of primitive $\hc$-invariant prime ideals in $\oqm2n$ is 
$(3^{n+1}-2^{n+1}+(-1)^{n+1}+2)/4$.  Cauchon \cite{c2} (see also \cite{laucomb}) enumerated the $\hc$-invariant prime ideals in $\oqmn$, giving a closed formula in terms of the Stirling numbers of the second kind.  In particular, the number of $\hc$-invariant prime ideals in $\oqm2n$ is $2\cdot 3^n - 2^n$.  
Surprisingly, these formulas show that the number of $\hc$-invariant prime ideals that are primitive is far from being negligible, and the proportion tends to $3/8$ as $n\rightarrow \infty$.  We also give a table of data obtained using Maple and some conjectures about $\hc$-invariant primitive ideals in $\oqmmn$.

Using our combinatorial criterion, one can show that the Gelfand-Kirillov dimension 
of every factor of $\oqmmn$ by an $\hc$-invariant primitive ideal is even. We next asked ourselves whether all primitive factors of $\oqmmn$ have even Gelfand-Kirillov dimension. It turns out that we are able to prove this result for a wide class of algebras---the so-called CGL extensions---by using both the $H$-stratification theory of Goodearl and Letzter, and the theory of deleting-derivations developed by Cauchon. Examples of CGL extensions are quantum affine spaces, the algebra of generic quantum matrices, the positive part $U_q^+(\mathfrak{g})$ of the quantised enveloping algebra of any semisimple complex Lie algebra, etc.  In particular, our result shows that the Gelfand-Kirillov dimension of the primitive
quotients of  the positive part $U_q^+(\mathfrak{g})$ of the quantised enveloping
algebra of any semisimple complex Lie algebra is always even, just as in the
classical case. Indeed, in the classical setting, it is a well-known Theorem of Dixmier that the primitive factors of enveloping algebras of
finite-dimensional complex nilpotent Lie algebras are isomorphic to Weyl algebras, and so have even Gelfand-Kirillov dimension. However, contrary to the classical situation, in the quantum case, primitive ideals are not always maximal and two primitive quotients (with the same even Gelfand-Kirillov dimension) are not always isomorphic; in the case where $\mathfrak{g}$ is of type $B_2$, there are three classes of primitive quotients of $U_q^+(B_2)$ of Gelfand-Kirillov dimension $2$ \cite{b2}.

The paper is organised as follows. In the first section, we recall the notion of CGL extension that was introduced in \cite{llr}. The advantage of these algebras is that one can use both the $H$-stratification theory of Goodearl and Letzter, and the deleting-derivations theory of Cauchon to study their prime and primitive spectra. After recalling, these two theories, we prove that every primitive factor of a (uniparameter) CGL extension has even Gelfand-Kirillov dimension. 

The second part of this paper is devoted to a particular (uniparameter) CGL extension: the algebra $\oqmmn$ of generic quantum matrices. We first prove our combinatorial criterion for an $\hc$-invariant prime ideal to be primitive. Then we use this criterion in order to obtain a formula for the total number of  primitive $\hc$-invariant ideals in $\oqm2n$. Finally, we give a table of data obtained using Maple and some conjectures about $\hc$-invariant primitive ideals in quantum matrices.

Throughout this paper, we use the following conventions. 
\\$\bullet$ If
$I$ is a finite set, $|I|$ denotes its cardinality.  
\\$\bullet$  $\gc a,b \dc := \{ i\in{\mathbb N} \mid a\leq i\leq b\}$. 
\\$\bullet$ $\mathbb{K}$
denotes a field and we set
$\mathbb{K}^*:=\mathbb{K}\setminus \{0\}$.  
\\$\bullet$ If $A$ is a $\C$-algebra, then $\spec(A)$ and  $\prim(A)$ denote respectively its prime and primitive spectra.

\section{Primitive ideals of CGL extensions.}
\label{SectionCGL}

In this section, we recall the notion of CGL extension that was introduced in \cite{llr}. 
Examples include various quantum algebras in the generic case such as quantum affine spaces, quantum matrices, 
positive part of quantised enveloping algebras of semisimple complex Lie algebras, etc. As we will see, the advantage of this class of algebras is that one can both use the stratification theory of Goodearl and Letzter and the theory of deleting-derivations of Cauchon in order to study their prime and primitive spectra. This will allow us to prove that every primitive factor of a (uniparameter) CGL extension has even Gelfand-Kirillov dimension.  

\subsection{$H$-stratification theory of Goodearl and Letzter, and CGL extensions.}
\label{sectionCGL}

Let $A$ denote a $\mathbb{K}$-algebra and $H$ be  a group acting on
 $A$ by $\mathbb{K}$-algebra automorphisms. A nonzero element $x$ of $A$ is an $H$-eigenvector of $A$ if 
$h\cdot x \in \mathbb{K}^{*}x$ for all $h \in H$. In this case, there
exists a character $\lambda$ of $H$ such that $h\cdot x=\lambda(h)x$ for
all $h \in H$, and $\lambda$ is called the $H$-eigenvalue of $x$. 

A two-sided ideal $I$ of $A$ is said to be $H$-invariant if $h\cdot I=I$ for all
$h \in H$.  An $H$-prime ideal of $A$ is a proper $H$-invariant ideal
$J$ of $A$ such that whenever $J$ contains the product of two
$H$-invariant ideals of $A$, $J$ contains at least one of them. We denote
by $H$-$\spec(A)$ the set of all $H$-prime ideals of $A$. Observe
that, if $P$ is a prime ideal of $A$, then
\begin{equation}
(P:H)\ := \ \bigcap_{h\in H} h\cdot P
\end{equation}
 is an $H$-prime ideal of
$A$. This observation allowed Goodearl and Letzter \cite{gl2}
(see also \cite{bgbook}), to construct a
partition of the prime spectrum of $A$ that is indexed by the
$H$-spectrum. Indeed, let $J$ be an $H$-prime ideal of $A$. We denote
by $\spec_J (A)$ the $H$-stratum associated  to $J$; that is,  
\begin{equation}
\spec_J (A)=\{ P \in \spec(A) \mbox{ $\mid$ } (P:H)=J \}.
\end{equation}
Then the $H$-strata of $\spec(A)$ form a partition of $\spec(A)$
\cite[Chapter II.2]{bgbook}; that
is:
\begin{equation}
 \spec(A)= \bigsqcup_{J \in H\mbox{-}\spec(A)}\spec_J(A).
 \end{equation}
This partition is the so-called $H$-stratification of $\spec(A)$. 
When the $H$-spectrum of $A$ is finite this
partition is a powerful tool in the study of the prime spectrum of $A$. 
In the generic case most quantum algebras have a finite $H$-spectrum (for a suitable action of
a torus on the algebra considered). We now move to the situation
where the $H$-spectrum is finite.

Throughout this paragraph $N$ denotes a positive integer and
 $R$ is an iterated Ore extension; that is,
\begin{equation}
R\ = \ \C[X_1][X_2;\sigma_2,\delta_2]\dots[X_N;\sigma_N,\delta_N],
\end{equation}
 where $\sigma_j$ is an automorphism of the $\C$-algebra $R_{j-1}:=\C[X_1][X_2;\sigma_2,\delta_2]\dots[X_{j-1};\sigma_{j-1},\delta_{j-1}]$
 and  $\delta_j$ is a  $\C$-linear $\sigma_j$-derivation of
 $R_{j-1}$ for all $j\in \{ 2,...,N \}$. Thus $R$ is a noetherian domain.  Henceforth, we
 assume that, in the terminology of \cite{llr}, $R$ is a CGL extension.

\begin{defn}[\cite{llr}]
{\em 
The iterated Ore extension $R$ is said to be a \emph{CGL extension} if 
\begin{enumerate}
\label{hypofond}
\item For all $j \in \gc 2,N \dc$, $\delta_j$ is locally nilpotent;
\item For all $j \in \gc 2,N \dc$, there exists  $q_j \in \C^*$ such
  that $\sigma_j \circ \delta_j = q_j \delta_j \circ \sigma_j$ 
and, for all $i \in \gc 1,j-1 \dc$, there exists $\lambda_{j,i} \in
  \C^*$ such that $\sigma_j(X_i)=\lambda_{j,i} X_i$;
\item None of the $q_j$ ($2 \leq j \leq N$) is a root of unity;
\item There exists a torus $H=(\C^*)^d$ that acts rationally by
  $\C$-automorphisms on $R$ such that:
\\$\bullet$ $X_1,\dots,X_N$ are $H$-eigenvectors;
\\$\bullet$ The set $\{ \lambda \in \C^* \mbox{ $\mid$ } (\exists h
  \in H)(h\cdot X_1=\lambda X_1)\}$ is infinite;
\\$\bullet$ For all $j \in \gc 2,N \dc$, there exists $h_j \in H$ such
  that $h_j\cdot X_i=\lambda_{j,i}X_i$ if $1 \leq i < j$ 
and $h_j\cdot X_j=q_j X_j$.
\end{enumerate}
}
\end{defn}

Some of our results will only be available in the ``uniparameter
case''.
\begin{defn}
\label{uniparameter}
{\em 
Let $R$ be a CGL extension. We say that $R$ is a \emph{uniparameter CGL extension} if 
there exist an antisymmetric matrix $(a_{i,j}) \in \mathcal{M}_N(\mathbb{Z})$ and $q \in \C^*$ 
not a root of unity  such that $\lambda_{j,i}=q^{a_{j,i}}$ for all 
$1 \leq j < i \leq N$.}
\end{defn}

The following result was proved by Goodearl and Letzter.

\begin{thm}\cite[Theorem II.5.12]{bgbook}
Every $H$-prime ideal of $R$ is completely prime, so that $H$-$\spec(R)$ coincides with the set of $H$-invariant completely prime ideals of $R$. Moreover
  there are at most $2^N$ $H$-prime ideals in $R$.
\end{thm}

As a corollary, the $H$-stratification breaks down the prime spectrum
of $R$ into a finite number of parts, the $H$-strata. The geometric
nature of the $H$-strata is well known: each $H$-stratum is
homeomorphic to the scheme of irreducible varieties of a $\C$-torus \cite[Theorems II.2.13 and II.6.4]{bgbook}. For completeness, we mention that the
$H$-stratification theory is a powerful tool to recognise primitive
ideals.

\begin{thm}\cite[Theorem II.8.4]{bgbook}
\label{theostratification}
The primitive ideals of $R$ are exactly the primes of $R$ that are
maximal in their $H$-strata.
\end{thm}

\subsection{A fundamental example: quantum affine spaces.}

Let $N$ be a positive integer and $\Lambda=\left( \Lambda_{i,j}
\right) \in \mathcal{M}_N (\C^*)$ a multiplicatively antisymmetric
matrix; that is, $\Lambda_{i,j} \Lambda_{j,i}=\Lambda_{i,i}=1$
 for all  $i,j \in \gc 1,N \dc$. The quantum affine space associated
 to $\Lambda$ is denoted by  $\mathcal{O}_{\Lambda}(\C^N)=
\C_{\Lambda}[T_1,\dots,T_N]$; this is the $\C$-algebra generated by
  $N$ indeterminates $T_1,\dots,T_N$ subject to the relations $T_j T_i
  = \Lambda_{j,i} T_i T_j$ for all $i,j \in \gc 1,N \dc$. It is
  well known that $\mathcal{O}_{\Lambda}(\C^N)$ is an iterated Ore
  extension that we can write:
$$\mathcal{O}_{\Lambda}(\C^N)=\C [T_1][T_2;\sigma_2]\dots[T_N;\sigma_N],$$
where $\sigma_j$ is the automorphism defined by
$\sigma_j(T_i)=\lambda_{j,i} T_i$ for all $1\leq i < j \leq
N$. Observe that the torus $H=\left( \C ^* \right) ^N$ acts by
automorphisms on $\mathcal{O}_{\Lambda}(\C^N)$ 
via:
$$(a_1,\dots,a_N)\cdot T_i=a_iT_i \mbox{ for all }i\in \gc 1,N \dc \mbox{ and }(a_1,\dots,a_N) \in H.$$
 Moreover, it is well known (see for instance \cite[Corollary 3.8]{llr}) that
 $\mathcal{O}_{\Lambda}(\C^N)$ is a CGL extension 
 with this action of $H$. Hence $\mathcal{O}_{\Lambda}(\C^N)$ has at
 most $2^N$  $H$-prime ideals and they are all completely prime.

The $H$-stratification of 
$\spec\left( \mathcal{O}_{\Lambda}(\C^N) \right)$ has been entirely
described by Brown and Goodearl when the group $\langle \lambda_{i,j}
\rangle$ is torsion free  \cite{bgtrans} and next by Goodearl and Letzter
in the general case \cite{gl1}. We now recall their results.

Let $W$ denote the set of subsets of $\gc 1,N \dc$. If $w \in W$,
then we denote by $K_w$ the (two-sided) ideal of
$\mathcal{O}_{\Lambda}(\C^N)$ generated by the indeterminates $T_i$
with $i\in w$. It is easy to check that $K_w$ is an $H$-invariant completely prime ideal of 
$\mathcal{O}_{\Lambda}(\C^N)$.  

\begin{prop}\cite[Proposition 2.11]{gl1}
\label{proprappelHstratificationespacesaffinesquantiques}
The following hold:
\begin{enumerate} 
\item The ideals $K_w$ with $w\in W$ are exactly the $H$-prime ideals
  of $\mathcal{O}_{\Lambda}(\C^N)$.  Hence there are exactly $2^N$
  $H$-prime ideals in that case;
\item For all  $w\in W$, the
  $H$-stratum associated to $K_w$ is given by
$$\spec_{K_w}\left( \mathcal{O}_{\Lambda}(\C^N) \right) =\left\{ P \in \spec \left(\mathcal{O}_{\Lambda}(\C^N) \right) \mbox{ $\mid$ } P \cap \{T_i \mbox{ $\mid$ }i\in \gc 1,N \dc \}= \{ T_i \mbox{ $\mid$ }i\in w\}\right\} .$$
\end{enumerate}
\end{prop}

\subsection{The canonical  partition of $\spec(R)$.}
\label{canonicalembedding}
In this paragraph, $R$ denotes a CGL extension as in Section \ref{sectionCGL}. We present the canonical partition of $\spec(R)$
that was constructed by Cauchon \cite{c1}. This partition gives new
insights to the $H$-stratification of $\spec(R)$.

In order to describe the prime spectrum of $R$, Cauchon \cite[Section 3.2]{c1} has
constructed an algorithm called the \emph{deleting-derivations algorithm}.  The reader is referred to \cite{c1,c2} for more details on this algorithm. 
One of the interests of this algorithm is that it has allowed Cauchon 
to rely the prime spectrum of a CGL extension to the prime spectrum of a certain quantum affine space. More precisely, let $\Lambda=\left( \mu_{i,j} \right) \in \mathcal{M}_N(\C^*)$ be the
multiplicatively antisymmetric matrix whose entries are defined as follows. 
$$\mu_{j,i}=\left\{ \begin{array}{ll}
\lambda_{j,i} & \mbox{ if }i<j \\
1 & \mbox{ if } i=j \\
\lambda_{j,i}^{-1} & \mbox{ if }i> j,
\end{array}\right.$$
where the  $\lambda_{j,i}$ with $i<j$ are coming from the CGL extension structure of $R$ (see Definition in Section \ref{sectionCGL}). Then we set $\overline{R}:=\mathbb{K}_\Lambda [T_1,\dots,T_N]=\mathcal{O}_{\Lambda}(\C^N)$.\\

Using his deleting-derivations algorithm, Cauchon has shown \cite[Section 4.4]{c1} that there exists 
an (explicit) embedding, called the canonical embedding, $\varphi:\spec(R) \longrightarrow \spec(\overline{R})$. 
This canonical embedding allows the construction of a partition of $\spec(R)$ as follows.
 
We keep the notation of the previous sections. In particular, $W$ still denotes the set of all subsets of $\gc 1,N \dc$. If $w \in W$, we
set 
$$  \spec_{w}\left(R \right)=\varphi^{-1} \left(  \spec_{K_w}\left( \overline{R} \right) 
\right).$$
Moreover, we denote by $W'$ the set of those $w \in W$ such that $
\spec_{w}\left(R \right) \neq \emptyset$. Then it follows from the work of Cauchon
\cite[Proposition 4.4.1]{c1} that 
$$\spec(R) =\bigsqcup_{w \in W'} \spec_w (R) \mbox{ and }\mid W' \mid \leq \mid W \mid = 2^N.$$

This partition is called the canonical partition of $\spec(R)$; this
  gives another way to understand the $H$-stratification since Cauchon
  has shown \cite[Th\'eor\`eme 5.5.2]{c1} that these two partitions
  coincide. As a consequence, he has given another description of
  the $H$-prime ideals of $R$.

\begin{prop}\cite[Lemme 5.5.8 and Th\'eor\`eme 5.5.2]{c1}
\label{descriHprem}
$ $
\begin{enumerate}
\item Let $w \in W'$. There exists a (unique) $H$-invariant (completely) prime ideal $J_w$ of
  $R$ such that $\varphi(J_w)=K_w$, where $K_w$ denotes the ideal of $\overline{R}$ generated by the $T_{i}$  with $i \in w$.
\item $ H \mbox{-}\spec(R) = \{J_w \mbox{ $\mid$ }w\in W'\}$. 
\item $\spec_{J_w}(R)=\spec_w(R)$ for all $w \in W'$.
\end{enumerate}
\end{prop}

Regarding the primitive ideals of $R$, one can use the canonical
embedding to characterise them. Indeed, let $P$ be a primitive ideal
of $R$. Assume that $P \in \spec_w(R)$ for some $w \in W'$. Then, it follows from Theorem \ref{theostratification} that 
$P$ is maximal in $ \spec_w(R)$. Now, recall from the work of Cauchon \cite[Th\'eor\`emes 5.1.1 and
  5.5.1]{c1} that the canonical
embedding induces an inclusion-preserving homeomorphism from $ \spec_w(R)$ onto
$\spec_{w}\left( \overline{R}
\right)=\spec_{K_w}(\overline{R})$. Hence $\varphi(P)$ is a maximal
ideal within $\spec_{w}\left( \overline{R}
\right)=\spec_{K_w}(\overline{R})$, and so we deduce from Theorem
\ref{theostratification} that $\varphi(P)$ is a primitive ideal of $\overline{R}$ that belongs
to   $\spec_{w}\left( \overline{R}
\right)=\spec_{K_w}(\overline{R})$. 
Also, similar arguments show that, if $P$ is a prime ideal of $R$ such
that $\varphi(P)$ is a primitive ideal of $\overline{R}$, then $P$ is
primitive. So, one can state the following result.

\begin{prop}
\label{primitivecanonicalembedding}
Let $P \in \spec(R)$ and assume that $P \in \spec_w(R)$ for some $w
\in W'$. Then $\varphi(P) \in \spec_{K_w}(\overline{R})$ and 
$P$ is primitive if and only if $\varphi(P)$ is primitive.
\end{prop}

This result was first obtained by Cauchon \cite[Th\'eor\`eme 5.5.1]{c1preprint}.

\subsection{Gelfand-Kirillov dimension of primitive quotients of a CGL extension.}
$ $

In this paragraph, $R$ still denotes a CGL extension. We start by recalling the notion of 
$\Tdeg$-stable algebra defined by Zhang \cite{z1}.

\begin{defn}
{\em
Let $A$ be a $\C$-algebra and $\mathcal{V}$ be the set of
finite-dimensional subspaces of $A$ that contain $1$.
\begin{enumerate}
\item Let $V \in \mathcal{V}$ and $n$ be a
  nonnegative integer. If $\{v_1,\dots,v_m\}$ is a basis of $V$, then
  we denote by $V^n$ the subspace of $A$ generated by the $n$-fold products of elements in $V$. (Here we use the
  convention $V^0=\C$.)
\item The \emph{Gelfand-Kirillov dimension} of $A$, denoted $\gkdim(A)$, is
  defined by:
$$\gkdim(A)=\sup_{V \in \mathcal{V}} \overline{\lim_{n \rightarrow \infty}} \frac{\log ({\rm dim}(V^n))}{\log (n)}.$$
\item The \emph{Gelfand-Kirillov transcendence degree} of $A$, denoted
  $\Tdeg(A)$, is defined by:
$$\Tdeg(A)=\sup_{V\in \mathcal{V}} \inf_{b} 
\overline{\lim_{n \rightarrow \infty}} \frac{\log ({\rm dim}((bV)^n))}{\log (n)},$$
where $b$ runs through the set of regular elements of $A$.
\item $A$ is  $\Tdeg$-\emph{stable} if the following hold:
\\$\bullet$ $\gkdim(A)=\Tdeg(A)$. 
\\$\bullet$ For every multiplicative system of regular elements $S$ of
$A$ that satisfies the Ore condition, we have:  $\Tdeg(S^{-1}A)=\Tdeg(A)$.
 \end{enumerate}}
\end{defn} 
$ $

Let $P \in \prim(R) \cap \spec_w(R)$ for some $w\in W'$. Then, it
follows from Proposition \ref{primitivecanonicalembedding} that $\varphi(P)$ is a primitive ideal of $\overline{R}$ that belongs to $\spec_{K_w}(R)$, where 
$$K_w=\langle T_i \mbox{ $\mid $ } i \in w \rangle .$$ 
Let $i \notin w$. We denote by $t_i$ the canonical image of $T_i$ in
the algebra $\overline{R}/K_w$. Also, we denote
by $B_w$ the subalgebra of $\fract \left( \overline{R}/K_w \right) $ defined by $B_w:= \C \langle t_i^{\pm 1} \mbox{ $\mid $ } i\notin w \rangle $. 
$B_w$ is the quantum torus associated to the quantum affine space $\overline{R}/K_w$. In other words, $B_w$ is a McConnell-Pettit algebra in $t_i$ with $i \notin w$ (see \cite{mcp}).

It follows from the work of Cauchon \cite[Th\'eor\`eme 5.4.1]{c1} that there exists a
multiplicative system of regular elements $\mathcal{F}$ of
$R/P$ that satisfies the Ore condition in
$R/P$, and such that
$$\left(R/P\right)
\mathcal{F}^{-1}=\left(\overline{R}/\varphi(P)\right)
\mathcal{E}^{-1} \simeq \frac{B_w}{\varphi(P)\mathcal{E}^{-1}},$$
where $\mathcal{E}$ denotes the canonical image of the multiplicative system of
$\overline{R}$ generated by the normal elements $T_i$ with $i \notin
w$. (Observe that $\varphi(P) \cap \mathcal{E} = \emptyset$ since
$\varphi(P) \in \spec_{w}(\overline{R})$.)

As $\varphi(P)$ is a primitive ideal of $\overline{R}$, we deduce from
\cite[Theorem 2.3]{gl1} that $\varphi(P)\mathcal{E}^{-1}$ is a
primitive ideal of the quantum torus $B_w$. As all the primitive
ideals of $B_w$ are maximal \cite[Corollary 1.5]{gl1}, $\varphi(P)\mathcal{E}^{-1}$ is a
maximal ideal of $B_w$ and
$$ \varphi(P)\mathcal{E}^{-1} = \left\langle
\varphi(P)\mathcal{E}^{-1} \cap Z(B_w) \right\rangle,$$
where $Z(B_w)$ denotes the centre of $B_w$. Also, it follows
from  \cite[Corollary 1.5]{gl1} that $ \varphi(P)\mathcal{E}^{-1} \cap Z(B_w)$ is a maximal
ideal of $Z(B_w)$. Recall from \cite[1.3]{gl1} that $Z(B_w)$ is a commutative
Laurent polynomial ring over $\C$.


We now assume that $R$ is a uniparameter CGL extension. In this case, it follows \cite[Proposition 2.3]{richard} that $B_w/\left(\varphi(P)\mathcal{E}^{-1}\right)$
is isomorphic to a simple quantum torus, and its GK dimension is an even integer.

As a quantum torus is $\Tdeg$-stable \cite[Proposition 7.2]{z1}, so is
$B_w/\left(\varphi(P)\mathcal{E}^{-1}\right)$. Moreover, as $\gkdim \left( B_w/\left(\varphi(P)\mathcal{E}^{-1}\right)
  \right)$ is an even integer, we see that  $\left(R/P\right)
  \mathcal{F}^{-1}$ is also  $\Tdeg$-stable of even Gelfand-Kirillov
dimension. As $R/P$
is a subalgebra of $\left(R/P\right) \mathcal{F}^{-1}$
such that $\fract \left( R/P \right)= \fract
  \left(\left( R/P\right) \mathcal{F}^{-1} \right)$, we deduce from
\cite[Proposition 3.5 (4)]{z1} the following results.

\begin{thm} 
\label{theoTdegstable}
Assume that $R$ is a uniparameter CGL extension and let $P$ be a primitive ideal of $R$. 
\begin{enumerate}
\item $\displaystyle{R/P}$ is $\Tdeg$-stable.
\item $\gkdim \left( R/P\right)$ is an even integer.
\end{enumerate}\label{thm: unipar}
\end{thm} 

This result can be applied to several quantum algebras. In particular,
it follows from \cite[Lemma 6.2.1]{c1} that it can be applied to the 
positive part $U_q^+(\mathfrak{g})$ of the quantised enveloping
algebra of any semisimple complex Lie algebra when the parameter $q \in \mathbb{K}^*$ is not 
a root of unity. As a result,
every primitive quotient of $U_q^+(\mathfrak{g})$ has even
Gelfand-Kirillov dimension. Roughly speaking, this is a quantum counterpart of the theorem of Dixmier that asserts that the primitive
factor algebras of the enveloping algebra $U(\mathfrak{n})$ of a
finite-dimensional complex nilpotent Lie algebra $\mathfrak{n}$ are
isomorphic to Weyl algebras. Indeed, our result shows that, as in the
classical case, the Gelfand-Kirillov dimension of a primitive
quotient of $U_q^+(\mathfrak{g})$ is always an even integer.  In the
quantum case however, primitive ideals are not always maximal, and two primitive quotients with the same Gelfand-Kirillov
dimension are not always isomorphic. Indeed, in the case where
$\mathfrak{g}$ is of type $B_2$, it turns out that there are three
classes of primitive quotients of $U_q^+(B_2)$ of Gelfand-Kirillov
dimension $2$ \cite{b2}. \\

\begin{remark}
 {\em The uniparameter hypothesis is needed in Theorem \ref{theoTdegstable}. Indeed, let ${\bf q}$ be any $3 \times 3$  multiplicatively antisymmetric matrix whose entries generate a free abelian group of rank 3 in $\mathbb{K}^*$. Then it follows from \cite[Proposition 1.3]{mcp} and \cite[Theorem 2.3]{gl1} that $(0)$ is a primitive ideal in the quantum affine space $\mathcal{O}_{\bf{q}} (\mathbb{K}^3)$. However the Gelfand-Kirillov dimension of $\mathcal{O}_{\bf{q}} (\mathbb{K}^3)$ is equal to 3, and so is not even.}
\end{remark}

\section{Primitive $\hc$-primes in quantum matrices.}

In this section, we study the primitive ideals of a particular CGL extension: the algebra of generic quantum matrices. In particular, we prove a combinatorial criterion for an $H$-prime ideal to be primitive in this algebra. Then we use this criterion to compute the number of primitive $H$-primes in the algebra of $2 \times n$ quantum matrices. The motivation to obtain a formula for the total number of primitive $H$-primes in the algebra of $m \times n$ quantum matrices comes from the fact that this number corresponds to the number of ``$H$-invariant points" in the ``variety of quantum matrices". We finish by giving some data and some conjectures for the number of primitive $H$-primes in $m \times n$ quantum matrices.

Throughout this section, \textbf{$q\in \mathbb{K}^*$ is not a root of unity}, and $m,n$ denote positive
integers.

\subsection{Quantum matrices as a CGL extension.}

We denote by $R=\oqmmn$ the standard quantisation of the ring of
regular functions on $m \times n$ matrices with entries in $\mathbb{K}$; it is
the $\mathbb{K}$-algebra generated by the $m \times n $ indeterminates
$Y_{\ia}$, $1 \leq i \leq m$ and $ 1 \leq \alpha \leq n$, subject to the
following relations:\\ \[
\begin{array}{ll}
Y_{i, \beta}Y_{i, \alpha}=q^{-1} Y_{i, \alpha}Y_{i ,\beta},
& (\alpha < \beta); \\
Y_{j, \alpha}Y_{i, \alpha}=q^{-1}Y_{i, \alpha}Y_{j, \alpha},
& (i<j); \\
Y_{j,\beta}Y_{i, \alpha}=Y_{i, \alpha}Y_{j,\beta},
& (i <j,  \alpha > \beta); \\
Y_{j,\beta}Y_{i, \alpha}=Y_{i, \alpha} Y_{j,\beta}-(q-q^{-1})Y_{i,\beta}Y_{j,\alpha},
& (i<j,  \alpha <\beta). 
\end{array}
\]

It is well known that $R$ can be presented as an iterated Ore extension over
$\mathbb{K}$, with the generators $Y_{\ia}$ adjoined in lexicographic order.
Thus the ring $R$ is a noetherian domain; we denote by $F$ its skew-field of
fractions. Moreover, since $q$ is not a root of unity, it follows from
\cite[Theorem 3.2]{gletpams} that all prime ideals of $R$ are completely
prime.

It is well known that the algebras $\oqmmn$ and $\oqmnm$ are
isomorphic. Hence, all the results that we will proved for $\oqm2n$ will also be valid for 
$\co_q(M_{n,2})$.

It is easy to check
that the group $\hc:=\left( \mathbb{K}^* \right)^{m+n}$ acts on $R$ by
$\mathbb{K}$-algebra automorphisms via:
$$(a_1,\dots,a_m,b_1,\dots,b_n).Y_{\ia} = a_i b_\alpha Y_{\ia} \quad {\rm
for~all} \quad \: (\ia)\in \gc 1,m \dc \times \gc 1,n \dc.$$ 
Moreover, as $q$ is not a root of unity, $R$ endowed with this action of $\hc$ is a CGL extension (see for instance \cite{llr}). This implies in particular that $\hc$-$\spec(R)$ is finite and that every $\hc$-prime is completely prime.

\subsection{$\hc$-primes and Cauchon diagrams.}
\label{sectionCauchondiagram}

As $R=\oqmmn$ is a CGL extension, one can apply the results of Section \ref{SectionCGL} to this algebra. In particular, using the theory of deleting-derivations, Cauchon has given a combinatorial description 
of $\hc$-$\spec(R)$. More precisely, in the case of the algebra $R=\oqmmn$, he has described the set $W'$ that appeared in Section \ref{canonicalembedding} as follows.

First, it follows from \cite[Section 2.2]{c2} that the quantum affine space $\overline{R}$ that appears in Section \ref{canonicalembedding} is in this case $\overline{R}=\mathbb{K}_{\Lambda}[T_{1,1},T_{1,2},\dots,T_{m,n}]$, 
where $\Lambda$ denotes the $mn \times mn$ matrix defined as follows. We set
$$A:=\left(
\begin{array}{ccccc}
 0 & 1 & 1 & \dots & 1 \\
-1 & 0 & 1 & \dots  & 1 \\
\vdots & \ddots &\ddots&\ddots &\vdots \\
-1 & \dots & -1 & 0 & 1  \\
 -1& \dots& \dots & -1 & 0 \\  
\end{array} 
\right) 
\in \mathcal{M}_{m}(\mathbb{Z}) \subseteq \mathcal{M}_{m}(\mathbb{C}), $$
and 
 $$B=(b_{k,l}):=\left( 
\begin{array}{ccccc}
 A & I_m & I_m & \dots & I_m \\
-I_m & A & I_m & \dots  & I_m \\
\vdots & \ddots &\ddots&\ddots &\vdots \\
-I_m & \dots & -I_m & A & I_m  \\
 -I_m& \dots& \dots & -I_m & A \\  
\end{array} 
\right)
\in \mathcal{M}_{mn}(\mathbb{C}), $$
where $I_m$ denotes the identity matrix of $\mathcal{M}_m$. Then $\Lambda $ is
the $mn \times mn$ matrix whose entries are defined by $\Lambda_{k,l}
=q^{b_{k,l}}$ for all $k,l \in \gc 1, mn \dc$.

We now recall the notion of Cauchon diagrams that first appears in \cite{c2}.  
\begin{defn}{\em
An $m\times n$ \emph{Cauchon diagram} $C$ is simply an $m\times n$ grid consisting of $mn$ squares in which certain squares are coloured black.  We require that the collection of black squares have the following property:
\vskip 2mm
\noindent If a square is black, then either every square strictly to its left is black or every square strictly above it is black.\\
We let $\Ca$ denote the collection of $m\times n$ Cauchon diagrams.}
\end{defn}
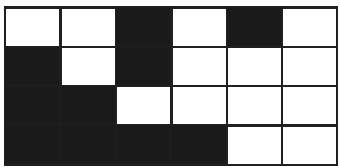
Figure 1: An example of a $4\times 6$ Cauchon diagram.
\vskip 2mm

\begin{tabular}{|p{0.30cm}|p{0.30cm}|p{0.30cm}|p{0.30cm}|p{0.30cm}|p{0.30cm}|}

\hline
 &  & \gray &  & \gray & 
 \\
\hline
\gray & & \gray &  & & 
 \\
\hline
\gray & \gray & & & & \\
\hline
\gray & \gray & \gray & \gray &  & \\
\hline
\end{tabular}
\vskip 2mm
Using the canonical embedding (see Section \ref{canonicalembedding}), Cauchon \cite{c2} produced a bijection between $\hc$-$\spec(\oqmmn)$ and the collection $\mathcal{C}_{m,n}$ of $m\times n$ Cauchon diagrams. Roughly speaking, with the notation of previous sections, the set $W'$ is the set of $m \times n$ Cauchon diagrams. Let us make 
this precise. If $C$ is a $m\times n$ Cauchon diagram, then we denote by $K_C$ the (completely) prime ideal 
of $\overline{R}$ generated by the indeterminates $T_{i,\alpha}$ such that the square in position $(i,\alpha)$ is a black square of $C$. Then, with $\varphi : \spec(R) \rightarrow \spec(\overline{R})$ denoting the canonical embedding, it follows from \cite[Corollaire 3.2.1]{c2} that there exists a unique $\hc$-invariant (completely) prime ideal $J_C$ of $R$ such that $\varphi (J_C)=K_C$; moreover there is no other $\hc$-prime in $\oqmmn$:
$$\hc \mbox{-} \spec (\oqmmn)=\{ J_C | C \in \mathcal{C}_{m,n} \}.$$

\begin{defn}
{\em
A Cauchon diagram $C$ is \emph{labeled} if each white square in $C$ is labeled with a positive integer such that:
\begin{enumerate}
\item{the labels are strictly increasing from left to right along rows;}
\item{if $i<j$ then the label of each white square in row $i$ is strictly less than the label of each white square in row $j$.}
\end{enumerate}
}
\end{defn}
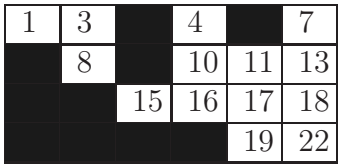
Figure 2: An example of a $4\times 6$ labeled Cauchon diagram.
\vskip 2mm

\begin{tabular}{|p{0.30cm}|p{0.30cm}|p{0.30cm}|p{0.30cm}|p{0.30cm}|p{0.30cm}|}

\hline
1 & 3 & \gray & 4 & \gray & 7
 \\
\hline
\gray & 8 & \gray & 10 & 11 & 13
 \\
\hline
\gray & \gray & 15 & 16 & 17 & 18 \\
\hline
\gray & \gray & \gray & \gray & 19 & 22 \\
\hline
\end{tabular}

\subsection{Perfect Matchings, Pfaffians and Primitivity.}

Our main tool in deducing when an $\hc$-prime ideal is primitive is to compute the \emph{Pfaffian} of a skew-symmetric matrix. We start with some background.  

\begin{notn}
 {\em Let $C$ be an $m\times n$ labeled Cauchon diagram with $d$ white squares and 
labels $\ell_1< \cdots < \ell_d$. 
\begin{enumerate}
 \item $A_C$ denotes the $d\times d$ skew-symmetric matrix whose $(i,j)$ 
entry is $+1$ if the square labeled $\ell_i$ is in the same column and strictly above the square labeled 
$\ell_j$ or is in the same row and strictly to the left of the square labeled $\ell_j$; its $(i,j)$ entry is $-1$ 
if the square labeled $\ell_i$ is in the same column and strictly below the square labeled $\ell_j$ or is in 
the same row and strictly to the right of the square labeled $\ell_j$; otherwise, the $(i,j)$ entry is $0$.
\item $G(C)$ denotes the directed graph whose vertices are the white squares of $C$ 
and in which we draw an edge from one white square to another if the first white square is either in the same row and strictly on the left of the second white square or the first white square is in the same column and strictly above the second white square.
\end{enumerate}
}
\end{notn}

Observe that $A_C$ is the skew adjacency matrix of the directed graph $G(C)$, and that both $A_C$ and $G(C)$ are independent of the set of labels which 
appear in $C$. Hence $A_C$ and $G(C)$ are defined for any Cauchon diagram.

\begin{defn}
 {\em Given a (labeled) Cauchon diagram $C$, the \emph{determinant of} $C$ is the element $\det (C)$ of $\mathbb{C}$ defined by:}
 $$\det (C) := \det(A_C).$$
\end{defn}

Before proving a criterion of primitivity for $J_C$ in terms of the Pfaffian of $A_C$ and perfect matchings, we first establish the following equivalent result.

\begin{thm}
 \label{thm: prim1} Let $P$ be an $\hc$-prime in $\oqmmn$.  Then $P$ is primitive if and only if the determinant of the Cauchon diagram corresponding to $P$ is nonzero.
\end{thm}
\begin{proof} Let $C$ be an $m\times n$ Cauchon diagram with $d$ white squares. We make $C$ into a labeled Cauchon diagram with labels $\ell_1< \cdots < \ell_d$. It follows from Proposition \ref{primitivecanonicalembedding} that $J_C$ is primitive if and only if 
$K_C$ is a primitive ideal of the quantum affine space $\overline{R}$, that is if and only if 
the ring $\frac{\overline{R}}{K_C}$ is primitive. Recall that  $\overline{R}=\mathbb{K}_{\Lambda}[T_{1,1},T_{1,2},\dots,T_{m,n}]$,
where $\Lambda$ denotes the $mn \times mn$ matrix whose entries are defined by $\Lambda_{k,l}
=q^{b_{k,l}}$ for all $k,l \in \gc 1, mn \dc$---the matrix $B$ has been defined in Section \ref{sectionCauchondiagram}. Let $\Lambda_C$ denote the multiplicatively antisymmetric $d \times d$ matrix whose entries are defined by $(\Lambda_C)_{i,j}=q^{(A_C)_{i,j}}$. 

As $K_C$ is the prime ideal generated by the indeterminates $T_{i,\alpha}$ such that the square in position $(i,\alpha)$ is a black square of $C$, the algebra $\frac{\overline{R}}{K_C}$ is isomorphic to 
the quantum affine space  $\mathbb{K}_{\Lambda_C}[t_1, \dots , t_d ]$ by an isomorphism that sends 
$T_{i,\alpha} +K_C $ to $t_k$ if the square of $C$ in position $(i,\alpha)$ is the white square labeled $\ell_k$, 
and $0$ otherwise.

Hence $J_C$ is primitive if and only if the quantum affine space $\mathbb{K}_{\Lambda_C}[t_1, \dots , t_d ]$ 
is primitive. To finish the proof, we use the same idea as in \cite[Corollary 1.3]{laulen}.

It follows from \cite[Theorem 2.3 and Corollary 1.5]{gl1} that 
the quantum affine space $\mathbb{K}_{\Lambda_C}[t_1, \dots , t_d ]$ 
is primitive if and only if the corresponding quantum torus 
$P(\Lambda_C):= \mathbb{K}_{\Lambda_C}[t_1, \dots , t_d ] \Sigma^{-1}$ is simple, 
where $\Sigma$ denotes the multiplicative system of 
$\mathbb{K}_{\Lambda_C}[t_1, \dots , t_d ]$ generated by the normal elements 
$t_1$, ..., $t_d$. Next, 
$\spec(P(\Lambda_C))$ is Zariski-homeomorphic via extension and contraction to
the prime spectrum of the centre $Z(P(\Lambda_C))$ of $P(\Lambda_C)$, by
\cite[Corollary 1.5]{gl1}. Further,
$Z(P(\Lambda_C))$ turns out to be a Laurent polynomial ring. To make 
this result precise, we need to introduce the following notation.

If $\underline{s}=(s_1,\dots,s_d) \in \mathbb{Z}^{d}$, then
we set $t^{\underline{s}}:= t_1^{s_1}\dots t_d^{s_d} \in P( \Lambda)$.
As in \cite{gl1}, we denote by $\sigma : \mathbb{Z}^{d} \times
\mathbb{Z}^{d} \rightarrow \mathbb{K}^*$ the antisymmetric bicharacter
defined by 
$$\sigma(\underline{s},\underline{t}):=\prod_{k,l=1}^d
(\Lambda_C)_{k,l}^{s_kt_l} \: \quad {\rm for~all}\quad  
\underline{s},\underline{t} \in
\mathbb{Z}^{d}.$$ 
Then it follows from \cite[1.3]{gl1} that the centre
$Z(P(\Lambda_C))$ of $P(\Lambda_C)$ is a Laurent polynomial ring over $\mathbb{K}$ in the variables
$(t^{\underline{b_1}})^{\pm 1},\dots,(t^{\underline{b_r}})^{\pm 1}$, where
$(\underline{b_1},\dots,\underline{b_r})$ is any basis of $S:=\{\underline{s}
\in \mathbb{Z}^{d} \mid \sigma(\underline{s},-)\equiv 1\}$. Since $q$ is not
a root of unity, easy computations show that $\underline{s} \in S$ if and only
if $A_C^t\underline{s}^t=0$. Hence the centre $Z(P(\Lambda_C))$ of $P(\Lambda_C)$
is a Laurent polynomial ring in $\dim_{\mathbb{C}} (\ker(A_C^t))$ indeterminates (here we
use the fact that $\dim_{\mathbb{Q}} (\ker(A_C^t))=\dim_{\mathbb{C}} (\ker(A_C^t))$). 
As a consequence, the quantum torus $P(\Lambda_C)$ is simple if and only if the matrix $A_C$ is invertible, that is, 
if and only if $\det (C) \neq 0$. To summarize, we have just proved that $J_C$ is primitive if and only if $\det (C) \neq 0$, as desired.\end{proof}

We finish this section by reformulating Theorem \ref{thm: prim1} in terms of Pfaffian and perfect matchings. Notice that the notion 
of perfect matchings of a directed graph or a skew-symmetric matrix is well known (see for instance \cite{Lov}). Roughly speaking, we define a perfect matching of a labeled Cauchon diagram as a perfect matching of the directed graph $G(C)$.

\begin{defn}
{\em Given a labeled Cauchon diagram $C$, we say that
$\pi=\{\{i_1,j_1\},\ldots ,\{i_m,j_m\}\}$ is a \emph{perfect matching} of $C$ if:
\begin{enumerate}
\item{$i_1,j_1,\ldots ,i_m,j_m$ are distinct;}
\item{$\{i_1,\ldots ,i_m,j_1,\ldots ,j_m\}$ is precisely the set of labels which appear in $C$;}
\item{$i_k<j_k$ for $1\le k\le m$;}
\item{for each $k$ the white square labeled $i_k$ is either in the same row or the same column as the white square labeled $j_k$.}
\end{enumerate}
We let $\PMC$ denote the collection of perfect matchings of $C$.}
\end{defn}

For example, for the Cauchon diagram in the Figure 2, we have the perfect matching
$\left\{\{1,4\}, \{3,8\}, \{7,13\}, \{10,16\}, \{11,17\},\{15,18\},\{19,22\}\right\}$. 

\begin{defn}{\em
Given a perfect matching
$\pi=\{\{i_1,j_1\},\ldots ,\{i_m,j_m\}\}$ of $C$ we call the sets $\{i_k,j_k\}$ for $k=1,\ldots ,m$ the \emph{edges} in $\pi$.
We say that an edge $\{i,j\}$ in $\pi$ is \emph{vertical} if the white squares labeled $i$ and $j$ are in the same column; otherwise we say that the edge is \emph{horizontal}.}
\end{defn}
Given a perfect matching $\pi$ of $C$, we define 
\begin{equation}
{\rm sgn}(\pi) \ := \ {\rm sgn}\left(\begin{array}{ccccccc} 1 & 2 & 3 & 4 & \cdots & 2m-1 & 2m\\
i_1 & j_1 & i_2 & j_2 & \cdots & i_m & j_m \end{array} \right).\end{equation}
We note that this definition of ${\rm sgn}(\pi)$ is independent of the order of the edges (see Lovasz \cite[p. 317]{Lov}).  It is vital, however, that $i_k<j_k$ for $1\le k\le m$.  
We then define 
\begin{equation} {\rm Pfaffian}(C) \ := \ \sum_{\pi\in \PMC} {\rm sgn}(\pi).\end{equation}
In particular, if $C$ has no perfect matchings, then $\Pf(C)=0$.

Observe that ${\rm Pfaffian}(C) $  is independent of the set of labels which appear in $C$, so that one can speak of the Pfaffian 
of any Cauchon diagram.

We are now able to establish the following criterion of primitivity for $J_C$. Even though this criterion is equivalent to the criterion given in Theorem \ref{thm: prim1}, this reformulation will be crucial in the following section.

\begin{thm} \label{thm: prim} Let $P$ be an $\hc$-prime in $\oqmmn$.  Then $P$ is primitive if and only if the Pfaffian of the Cauchon diagram corresponding to $P$ is nonzero.
\end{thm}
\begin{proof} Let $C$ be an $m\times n$ Cauchon diagram. It follows from Theorem \ref{thm: prim1} that $J_C$ is primitive if and only if 
the determinant of $A_C$ is nonzero. Since the determinant of $A_C$ is the square of the Pfaffian of $C$ \cite[Lemma 8.2.2]{Lov}, $J_C$ is primitive if and only if the Pfaffian of $C$ is nonzero, as claimed. \end{proof}

To compute the sign of a permutation, we use \emph{inversions}.  
\begin{defn} {\em Let ${\bf x}=(i_1,i_2,\ldots ,i_n)$ be a finite sequence of real numbers.  We define ${\rm inv}({\bf x})$ to be 
$\#\{ (j,k)~|~j<k, i_j>i_k\}$.  Given another finite real sequence ${\bf y}=(j_1,\ldots ,j_m)$, we define
${\rm inv}({\bf x}|{\bf y})=\#\{(k,\ell)\in \gc 1,n \dc \times \gc 1,m \dc~|~j_k<i_{\ell}\}$.}
\end{defn}
The key fact we need is that if $\sigma$ is a permutation in $S_n$, then
\begin{equation} {\rm sgn}(\sigma) \ = \ (-1)^{{\rm inv}(\sigma(1),\ldots ,\sigma(n))}.
\end{equation}

\subsection{Enumeration of $H$-primitive ideals in $\oqm2n$.}
 
In this section, we give a formula for the number of primitive $\hc$-prime ideals in the ring of $2\times n$ quantum matrices. We begin with some notation.
\begin{notn} {\em Given a statement $S$, we take $\chi(S)$ to be $1$ if $S$ is true and to be $0$ if $S$ is not true.}
\end{notn} 
 
We now compute the Pfaffian of a $2\times n$ Cauchon diagram.  To do this, we first must find the Pfaffian of a $1\times n$ Cauchon diagram. 
 \begin{lem} Let $C\in \mC_{1,n}$ be a $1\times n$ Cauchon diagram.  Then \[
 \Pf(C) \ = \ \left\{ \begin{array}{cl} 1 & {\rm if ~the~ number ~of~ white~ squares~ in~} C{\rm ~ is~ even;}\\
 0 & {\rm  otherwise.} \end{array}\right. \]
 \label{lem: 1cauch}
\end{lem}
\begin{proof} If the number of white squares in $C$ is odd, then $C$ has no perfect matchings and hence its Pfaffian is zero.  Thus we may assume that the number of white squares is an even integer $2m$ and the white squares are labeled from $1$ to $2m$ from left to right.  We prove that the Pfaffian is $1$ when $C$ has $2m$ white squares by induction on $m$.  When $m=1$, there is only one perfect matching and its sign is $1$.  Thus we obtain the result in this case.  We note that any perfect matching of $C$ is going to contain $\{1,i\}$ for some $i$.  Thus $\pi=\{1,i\}\cup \pi'$, where $\pi'$ is a perfect matching of the Cauchon diagram $C_i$ obtained by taking $C$ and colouring the white squares labeled $1$ and $i$ black.
Write $\pi'=\{\{i_1,j_1\},\ldots ,\{i_{m-1},j_{m-1}\}\}$ and let 
${\bf x}=(1,i,i_1,j_1,\ldots ,i_{m-1},j_{m-1})$ and let
${\bf x}'=(i_1,j_1,\ldots ,i_{m-1},j_{m-1})$.
Then ${\rm inv}({\bf x})={\rm inv}({\bf x}')+ (i-2)$.
Hence ${\rm sgn}(\pi)={\rm sgn}(\pi')(-1)^{i-2}$.
Since there is a bijective correspondence between perfect matchings of $C$ which contain $\{1,i\}$ and perfect matchings of $C_i$ we see that
$$\sum_{\{\pi\in \PMC~|~\{1,i\}\in\pi\}} {\rm sgn}(\pi) = 
\sum_{\pi'\in \PMCi} (-1)^{i-2}{\rm sgn}(\pi')=(-1)^{i-2}$$ by the inductive hypothesis.  Hence \begin{eqnarray*}
\Pf(C) &=& \sum_{\pi\in \PMC} {\rm sgn}(\pi) \\
&=& \sum_{i=2}^{2m} \sum_{\{\pi\in \PMC~|~\{1,i\}\in\pi\}} {\rm sgn}(\pi) \\
&=& \sum_{i=2}^{2m} (-1)^{i-2}\\
&= & 1. \end{eqnarray*}
 \end{proof}
\vskip 1mm
\noindent  In particular, we see that a $1\times n$ Cauchon diagram corresponds to a primitive ideal if and only if the number of white squares is even.  This is a special case of \cite[Theorem 1.6]{laulen}, but we need the value of the Pfaffian to study the $2\times n$ case.

\begin{notn}{\em 
 We let 
 $\Ca'$ denote the collection of $m\times n$ Cauchon diagrams which do not have any columns which consist entirely of black squares.  }
 \end{notn}
 We note that if $C\in \mC_{m,n}$ has exactly $d$ columns consisting entirely of black squares, then if we remove these $d$ columns we obtain an element of $\mC'_{m,n-d}$.  Hence we obtain the relation
 \begin{equation} 
 \label{eq: C}
 |\Ca| \ = \  \sum_{i=0}^n {n\choose i}|\mC'_{m,n-i}|, \end{equation}
 where we take $|\mC_{m,0}|=|\mC_{m,0}'|=1$.

 We begin by enumerating the elements of $\mC_{2,n}'$ which correspond to primitive $\hc$-primes in $\oqm2n$.  Let $C\in \mC_{2,n}'$.  Then the second row of $C$ has a certain number of contiguous black squares beginning at the lower left square.  If the second row does not consist entirely of black squares, then there is some smallest $i\ge 1$ for which the $(2,i)$ square of $C$ is white.  We note that the $(2,j)$ square must also be white for $i\le j\le n$, since otherwise we would necessarily have a column consisting entirely of black squares by the conditions defining a Cauchon diagram.    Since $C$ has no columns consisting entirely of black squares, the $(1,j)$ square is white for $1\le j<i$.  The remaining $n-i$ squares in positions $(1,j)$ for $i\le j\le n$ can be coloured either white or black and the result will still be an element of $\mC_{2,n}'$.  
Hence
\begin{equation} |\mC_{2,n}'| \ = \ \sum_{i=0}^n 2^{n-i} \ = \  2^{n+1}-1.\end{equation}
To enumerate the primitive $\mathcal{H}$-primes of $\mathcal{O}_q(M_{2,n})$, we need to introduce the following terminology.
\begin{notn}{\em
Given an element $C$ of $\mC_{2,n}'$, we use the following notation:
\begin{enumerate}
\item $p(C)$ denotes the largest $i$ such that the $(2,i)$ square of $C$ is black. 
\item  ${\rm Vert}(C)$ denotes the set of $j\in \{p(C)+1,\ldots ,n\}$ such that the $(1,j)$ square of $C$ is white.  (We use the name ${\rm Vert}(C)$, because this set consists of precisely the set of $j$ such that the $j$'th column of $C$ is completely white and hence it is only in these columns where a vertical edge can occur in a perfect matching of $C$.)
\item  Given a perfect matching $\pi$ of $C$ we let ${\rm Vert}(\pi)$ denote the set of $j\in {\rm Vert}(C)$ such that there is a vertical edge in $\pi$ connecting the two white squares in the $j$'th column.  
\item Given a subset $T\subseteq \{1,2,\ldots ,n\}$ we let ${\rm sum}_C(T)$ denote the sum of the labels in all white squares in the columns indexed by the elements of $T$. 
\end{enumerate}}
\end{notn}
 For example, if we use the Cauchon diagram $C$ in Figure 3 below, then $p(C)=1$, ${\rm Vert}(C)=\{2,4,6\}$, ${\rm sum}_C(\{2,3\})=(2+5)+6=13$.
 \begin{lem} Let $C\in \mC_{2,n}'$ be a labeled Cauchon diagram with $m$ white squares in the first row and $m'$ white squares in the second row with labels $\{1,2,\ldots ,m+m'\}$, and let $T$ be a subset of ${\rm Vert}(C)$.  Then
\[ \sum_{\{\pi\in \PMC~:~{\rm Vert}(\pi)=T\}} {\rm sgn}(\pi) =\left\{  
 \begin{array}{cl} (-1)^{{|T|+1\choose 2}+{\rm sum}_C(T)} & {\rm if}~m\equiv m'\equiv |T|~(\bmod 2); \\
0 & otherwise. \end{array}
\right.
 \]
 \label{lem: decomposition}
\end{lem}
\begin{proof} Let $C_1$ denote the first row of $C$ except that all squares in position $(1,j)$ with $j\in T$ are now coloured black and their labels are removed.  Similarly, let $C_2$ denote the second row of $C$ with the squares in positions $(2,j)$ with $j\in T$ coloured black and their labels removed (see Figure 3).  
\vskip 4mm
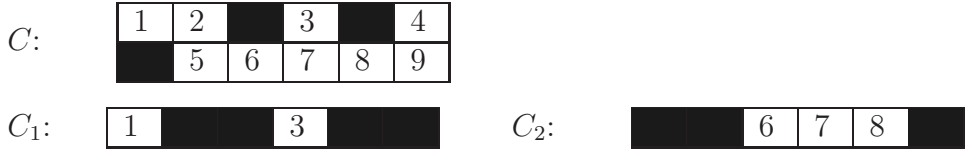
\noindent Figure 3: the decomposition of a $2\times 6$ labeled Cauchon diagram into rows with $T=\{2,6\}$.
\vskip 5mm
$C$:~~~~~~
\begin{tabular}{|p{0.30cm}|p{0.30cm}|p{0.30cm}|p{0.30cm}|p{0.30cm}|p{0.30cm}|}

\hline
1 & 2 & \gray & 3 & \gray & 4
 \\
\hline
\gray & 5 & 6& 7 & 8 & 9 \\
\hline
\end{tabular}
\vskip 3mm
$C_1$:~~~~
\begin{tabular}{|p{0.30cm}|p{0.30cm}|p{0.30cm}|p{0.30cm}|p{0.30cm}|p{0.30cm}|}
\hline
1 & \gray & \gray & 3 & \gray & \gray
 \\
\hline
\end{tabular}
~~~~~~$C_2$:~~~~~~
\begin{tabular}{|p{0.30cm}|p{0.30cm}|p{0.30cm}|p{0.30cm}|p{0.30cm}|p{0.30cm}|}
\hline
\gray & \gray & 6 & 7 & 8 & \gray
 \\
\hline

\end{tabular}

\vskip 3mm
\noindent 
Then the perfect matchings $\pi$ of $C$ with $\Vpi=T$ are in one-to-one correspondence with ordered pairs
$(\pi_1,\pi_2)$ in which $\pi_j$ is a perfect matching of $C_j$ for $j=1,2$.  Let $\pi_j$ be a perfect matching of $C_j$ for $j=1,2$.  It is therefore no loss of generality to assume that $C_1$ and $C_2$ both have an even number of squares.   We let $t=|T|$.  
Write \[\pi_1 = \{ \{a_1,b_1\}, \{a_2,b_{2}\},\ldots , \{a_r, b_{r}\} \}\]
 and
\[\pi_2 = \{ \{c_1,d_1\}, \{c_2,d_{2}\},\ldots , \{c_s, d_{s}\} \}.\]
We write \[\rho =\{ \{e_1,f_{1}\}, \{e_2,f_{2}\},\ldots , \{e_t, f_{t}\} \},\]
 where 
$e_1<e_2<\cdots <e_t$ are the labels of the white squares which appear in the positions
$\{(1,j)~|~j\in T\}$ and
$f_1<f_{2}<\cdots <f_{t}$ are the labels which appear in the positions $\{(2,j)~|~j\in T\}$.  We note that $\rho$ is precisely the vertical edges in the perfect matching $\pi=\pi_1\cup \pi_2\cup \rho$ of $C$.   
Let
\[ {\bf x}_1=(a_1,b_1,\ldots ,a_r,b_r), ~~~{\bf x}_2 = (c_1,d_1,\ldots ,c_s,d_s),~~~  {\bf x}_3=(e_1,f_{1},\ldots ,e_t,f_{t}). \]
Finally, let 
\[ {\bf x} = {\bf x}_1{\bf x}_3{\bf x}_2 = (a_1,b_1,\ldots ,a_r,b_r,e_1,f_{1},\ldots ,e_t,f_{t},c_1,d_1,\ldots ,c_s,d_s). \]
Note that \begin{equation}\label{eq: sgn} {\rm sgn}(\pi) = (-1)^{{\rm inv}({\bf x})}\qquad {\rm sgn}(\pi_i)=(-1)^{{\rm inv}(x_i)}\qquad {\rm for~}i=1,2.\end{equation}
We have
$${\rm inv}({\bf x}) = {\rm inv}({\bf x}_1)+{\rm inv}({\bf x}_2) + {\rm inv}({\bf x}_3) + {\rm inv}({\bf x}_1 | {\bf x}_2) + {\rm inv}({\bf x}_1 | {\bf x}_3) + {\rm inv}({\bf x}_3|{\bf x}_2).$$
Since
\[ e_1 < e_2 < \cdots < e_t < f_{1} < \cdots < f_{t}\] we have
$${\rm inv}({\bf x}_3)={\rm inv}(e_1,f_{1},e_2,f_{2},\ldots ,e_t,f_{t}) = (t-1)+\cdots +1={t\choose 2}.$$
Then $t+2r=m$ and $\{e_1,\ldots ,e_t,a_1,\ldots ,a_r,b_1,\ldots , b_r\} = \{1,2,\ldots ,m\}$. Also
$t+2s=m'$ and $\{f_{1},\ldots ,f_{t},c_1\ldots ,c_s,d_1,\ldots ,d_s\}=\{m+1,m+2,\ldots ,m+m'\}$. 
Notice that the $f_{1},\ldots ,f_{d}$ are greater than the labels appearing in $C_1$ and hence
$${\rm inv}({\bf x}_1|{\bf x}_3) =  \#\{(k,\ell)~|~e_k<a_{\ell}\} + \#\{(k,\ell)~|~e_k<b_{\ell}\}.$$
Since $\{a_1,\ldots ,a_r,b_1,\ldots ,b_r\}=\{1,2,\ldots ,m\}\setminus \{e_1,\ldots ,e_t\}$, we see that for each $k$,

$ \#\{\ell~|~e_k<a_{\ell}\} + \#\{\ell~|~e_k<b_{\ell}\}\ = \ \#\{e_k+1,\ldots ,m\}\setminus \{e_{k+1},\ldots ,e_t\}=m-e_k-(t-k)$.
Thus
\begin{eqnarray*} {\rm inv}({\bf x}_1|{\bf x}_3) &=& \sum_{k=1}^t (m-e_k+t-k) \\ &=&  mt - (e_1+\cdots + e_t) + {t\choose 2}.\end{eqnarray*}
To compute
${\rm inv}({\bf x}_3|{\bf x}_2)$, note that $e_1,\ldots ,e_t$ are all less than $c_1,d_1,\ldots ,c_s,d_s$ and hence
\begin{eqnarray*}
{\rm inv}({\bf x}_3|{\bf x}_2) 
&=& \sum_{k=1}^t \#\{\ell~|~c_{\ell}<f_{k}\}+\#\{\ell~|~d_{\ell}<f_{k}\}\\
 &=& \sum_{k=1}^t \#\{m+1,m+2,\ldots ,f_{k}\}\setminus \{f_{1},\ldots ,f_{{k}}\} \\
&=& \sum_{k=1}^t (f_k-m-k) \\
&=& -{t+1\choose 2}-mt+\sum_{k=1}^t f_{k}.\end{eqnarray*}
Note that ${\rm inv}({\bf x}_1|{\bf x}_2)=0$ and thus
\begin{eqnarray*}{\rm inv}({\bf x}) &=& {\rm inv}({\bf x}_1)+{\rm inv}({\bf x}_2) + {t\choose 2}  - (e_1+\cdots + e_t) + {t\choose 2}-{t+1\choose 2}+\sum_{k=1}^t f_{k}\\
&=& {\rm inv}({\bf x}_1)+{\rm inv}({\bf x}_2)+{t\choose 2}-t+{\rm sum}_C(T) -2 (e_1+\cdots + e_t) .\end{eqnarray*}
Equation (\ref{eq: sgn}) now gives
$${\rm sgn}(\pi)={\rm sgn}(\pi_1){\rm sgn}(\pi_2)(-1)^{{t+1\choose 2}+{\rm sum}_C(T)}.$$
It follows that 
\begin{eqnarray*} & & \sum_{\{\pi\in \PMC~:~\VC=T\}} {\rm sgn}(\pi)\\ &=&
\sum_{\pi_1\in {\mathcal PM}(C_1)}\sum_{\pi_2\in \mathcal{PM}(C_2)} 
{\rm sgn}(\pi_1){\rm sgn}(\pi_2)(-1)^{{t+1\choose 2}+{\rm sum}_C(T)}\\
&=& (-1)^{{t+1\choose 2}+{\rm sum}_C(T)}
\sum_{\pi_1\in {\mathcal PM}(C_1)} {\rm sgn}(\pi_1)
\sum_{\pi_2\in {\mathcal PM}(C_2)} {\rm sgn}(\pi_2) \\
&=& (-1)^{{t+1\choose 2}+{\rm sum}_C(T)}
\end{eqnarray*}
where the last step follows from Lemma \ref{lem: 1cauch}. 
\end{proof}

\begin{thm}  Let $C\in \mathcal{C}'_{2,n}$ be a Cauchon diagram with $m$ white squares in the first row and $m'$ white squares in the second row with labels $\{1,2,\ldots ,m+m'\}$, and let $S={\rm Vert}(C)$. Then ${\rm Pfaffian}(C)\not = 0$ if and only if $m\equiv m'\bmod 2$ and 
$$|S|-2{\rm sum}_C(S)\not\equiv 2m+2\bmod 4.$$
\label{thm: criterion}
\end{thm}
\begin{proof}  Let $S_0$ consist of the elements $j\in S$ with ${\rm sum}_C(\{j\})$ even and let $S_1$ consist of the elements $j\in S$ with ${\rm sum}_C(\{j\})$ odd.  
Notice that if $T\subseteq S$, then ${\rm sum}_C(T)\equiv |T\cap S_1|\bmod 2$.  Hence by Lemma \ref{lem: decomposition}
\begin{eqnarray*} \Pf(C) &=&
\sum_{\pi\in \PMC} {\rm sgn}(\pi) \\
&=& \sum_{T\subseteq S}~\sum_{\substack{\pi\in \PMC \\ \Vpi=T}} {\rm sgn}(\pi) \\
&=&\sum_{T\subseteq S} (-1)^{{|T|+1\choose 2}+{\rm sum}_C(T)}\chi(m\equiv m'\equiv |T|\bmod 2) \\
 &=& \sum_{T_0\subseteq S_0}\sum_{T_1\subseteq S_1} (-1)^{{|T_0|+|T_1|+1\choose 2}+|T_1|}
 \chi(m\equiv m'\equiv |T_0|+|T_1|\bmod 2) \\
 &=& \sum_{a=0}^{|S_0|}\sum_{b=0}^{|S_1|}{|S_0|\choose a}{|S_1|\choose b} (-1)^{{a+b+1\choose 2}+b}
 \chi(m\equiv m'\equiv a+b\bmod 2).\end{eqnarray*}
 At this point, we divide the evaluation of this sum into three cases.
 \vskip 1mm
 \noindent {\bf CASE I:} $m\not \equiv m'\bmod 2$.
 \vskip 1mm
 \noindent In this case, $ \chi(m\equiv m'\bmod 2)=0$ and thus $\Pf(C)=0$.\vskip 1mm
 \noindent {\bf CASE II:} $m$ and $m'$ are both odd.
 \vskip 1mm
 \noindent In this case, 
  \[ {\rm Pfaffian}(C) = \sum_{a=0}^{|S_0|}\sum_{b=0}^{|S_1|}{|S_0|\choose a}{|S_1|\choose b}(-1)^{{a+b+1\choose 2}+b}\chi(a+b\equiv 1\bmod 2).\]
 Note that if $a+b$ is odd, then ${a+b+1\choose 2}\equiv (a+b+1)/2\bmod 2$ and hence
 \begin{eqnarray*}
 \Pf(C) &=& \sum_{a=0}^{|S_0|}\sum_{b=0}^{|S_1|}{|S_0|\choose a}{|S_1|\choose b} (-1)^{(a+3b+1)/2}\chi(a+b\equiv 1\bmod 2) \\
 &=& {\rm Re}\Bigg( i \sum_{a=0}^{|S_0|}\sum_{b=0}^{|S_1|}{|S_0|\choose a}{|S_1|\choose b} i^a i^{3b}\Bigg) \\
 &=&  {\rm Re}\Bigg( i (1+i)^{|S_0|}(1-i)^{|S_1|}\Bigg) \\
 &=& \sqrt{2}^{|S_0|+|S_1|}{\rm Re}(i\exp(\frac{\pi i(|S_0|-|S_1|)}{4}))\\
 &=& \sqrt{2}^{|S_0|+|S_1|}\sin(\frac{\pi(|S_1|-|S_0|)}{4}).
 \end{eqnarray*}
Thus we see that the Pfaffian is $0$ in this case if and only if $|S_0|\equiv |S_1|\bmod 4$.  Notice that
$|S_0|-|S_1| = |S|-2|S_1|$.  Moreover, $|S_1|\equiv {\rm sum}_C(S)\bmod 2$ and hence the Pfaffian is $0$ if and only if 
$$|S|-2{\rm sum}_C(S)\equiv 0 \equiv 2m+2\bmod 4.$$
 \vskip 1mm
 \noindent {\bf CASE III:} $m$ and $m'$ are both even.
 \vskip 1mm
 \noindent The argument here is similar to the argument in Case II.  We now use the fact that if 
 $a+b\equiv 0 \bmod 2$, then $(-1)^{a+b+1\choose 2}= (-1)^{(a+b)/2}$.
 Hence
  \begin{eqnarray*}
 \Pf(C)  &=& {\rm Re}\Bigg( \sum_{a=0}^{|S_0|}\sum_{b=0}^{|S_1|}{|S_0|\choose a}{|S_1|\choose b} i^a i^{3b}\Bigg) \\
 &=& {\rm Re}\Bigg( (1+i)^{|S_0|}(1-i)^{|S_1|}\Bigg) \\
&=& \sqrt{2}^{|S_0|+|S_1|}\cos(\frac{\pi(|S_1|-|S_0|)}{4}) \end{eqnarray*}
Thus the Pfaffian is $0$ in this case if and only $|S_1|-|S_0|\equiv 2\bmod 4$. 
Again, we have 
$-|S_1|+|S_0|\equiv |S|-2{\rm sum}_C(S)$ and hence the Pfaffian of $C$ is $0$ if and only if
$$|S|-2{\rm sum}_C(S)\equiv 2 \equiv 2m+2\bmod 4.$$
\vskip 1mm
Thus we see that in each case we obtain the desired result. 
\end{proof}
\begin{thm} Let $n$ be a positive integer.  Then the number of primitive $\hc$-prime ideals in the ring $\oqm2n$ is $$\frac{3^{n+1}-2^{n+1}+(-1)^{n+1}+2}{4}.$$
\end{thm}
\begin{proof} For $a,b\in \{0,1\}^2$, we let $f_{a,b}(n)$ denote the number of Cauchon diagrams $C$ in $\mathcal{C}'_{2,n}$ with nonzero Pfaffian and with $p(C)\equiv a\bmod 2$ and $|{\rm Vert}(C)|\equiv b\bmod 2$. Then there are $p(C)+|{\rm Vert}(C)|\equiv b+a\bmod 2$ white squares in the first row of $C$ and $n-p(C)\equiv n-a \bmod 2$ white squares in the second row of $C$.  
We look at several cases. \vskip 1mm
\noindent {\bf CASE I:} $n+b$ is odd.  \vskip 2mm
In this case, the total number of white squares is odd and hence the Pfaffian is always zero.  Thus 
$f_{a,b}(n)=0$ in this case.
\vskip 1mm
\noindent {\bf CASE II:} $n$ and $b$ are odd ($b=1$).
\vskip 1mm
Notice that if $n+b$ is even and $b$ is odd, then
Theorem \ref{thm: criterion} gives automatically that we have nonzero Pfaffian since
in this case $|{\rm Vert}(C)|-2{\rm sum}_C ({\rm Vert}(C))-(2m+2)\equiv b\equiv 1\bmod 2$.  
Hence 
\begin{eqnarray*} f_{a,1}(n)&=& \sum_{p(C)=0}^n \sum_{i=0}^{n-p(C)} {n-p(C)\choose i}\chi(i\equiv 1\bmod 2)\\
&=& \sum_{j=0}^{n-1} 2^{n-j-1}\\
&=& 2^n-1,\end{eqnarray*}
 if $b$ and $n$ are both odd.   
\vskip 1mm
\noindent {\bf CASE III:} $n$ is even, $b$ is even, $a$ is odd ($(a,b)=(1,0)$).
\vskip 2mm
Let us start by looking at Cauchon diagrams in $\mathcal{C}'_{2,n}$ with $p(C)=a'\equiv 1\bmod 2$ and
$|{\rm Vert}(C)|=b'\equiv 0\bmod 2$. Such a Cauchon diagram is labeled with the trivial label: the labels on the first row are $1, \dots, a'+b'$ and the labels in the second row of $C$ are $\{a'+b'+1,\ldots ,b'+n\}$.  Let $J$ be a subset of this interval of even size $b'$. Let $C_J$ denote the Cauchon diagram in $\mathcal{C}'_{2,n}$ with $p(C)=a'$ and $S_J={\rm Vert}(C_J)$ consisting of all $j$ such that the $(2,j)$ entry of $C_J$ is a white square with label in $J$.  The labels of the squares in positions $(1,j)$ with $j\in {\rm Vert}(C_J)$ are just
$a'+1,\ldots ,a'+b'$. Then $${\rm sum}_{C_J}(S_J)=(a'+1)+\cdots +(a'+b')+\sum_{j\in J} j.$$
Since $b'$ is even, $$(a'+1)+\cdots + (a'+b')\equiv b'(b'+1)/2 \equiv b'/2\bmod 2.$$  By Theorem \ref{thm: criterion}, a necessary and sufficient condition for the Pfaffian to be nonzero in this case is
$$-b'+2{\rm sum_{C_J}}(S_J)-2(a'+b')-2 \not \equiv 0 \bmod 4.$$
Note that $-b'+2{\rm sum_{C_J}}(S_J)\equiv 2\sum_{j\in J} j\bmod 4$ and since $b'$ is even and $a'$ is odd, we see that
the Pfaffian is nonzero in this case if and only if
$$\#\{j\in J~|~j\equiv 1\bmod 2\}$$ is odd.  Note that $\{b'+a'+1,\ldots , n+b'\}$ is a set with
$n-a'$ elements, $(n-a'+1)/2$ are even and $(n-a'-1)/2$ are odd.  The number of ways of choosing a set $J$ of size $b'$ with $\#\{j\in J~|~j\equiv 1\bmod 2\}$ odd is then
$$\sum_{i=0}^{b'} {(n-a'-1)/2\choose i}{(n-a'+1)/2\choose b'-i}\chi(i\equiv 1\bmod 2).$$
Thus 
\begin{eqnarray*}
f_{1,0}(n) &=& \sum_{a'+b'\le n} \sum_{i=0}^{b'} {(n-a'-1)/2\choose i}{(n-a'+1)/2\choose b'-i}\chi(i-1 \equiv b'\equiv a'-1\equiv 0\bmod 2) \\
&=& \sum_{j=1}^{n/2} \sum_{i=0}^{n/2-j}\sum_{k=0}^{n/2-j+1} {n/2 - j\choose i}{n/2-j+1\choose k}\chi(i\equiv k\equiv 1\bmod 2) \\
&=& \sum_{j=1}^{n/2-1} 2^{n/2-j-1}2^{n/2-j} \\
&=& \sum_{j=1}^{n/2-1} 2^{n-2j-1} \\
&=& (2+2^3+\cdots + 2^{n-3}) \end{eqnarray*}
\vskip 1mm
\noindent {\bf CASE IV:} $n$, $a$ and $b$ are even ($(a,b)=(0,0)$).
\vskip 1mm
This case is treated using similar arguments than in case III. We keep the notation of case III. In particular, the labels in the second row of $C$ are just $\{b'+a'+1,\ldots ,n+b'\}$. Again, we must select a subset $J$ of size $b'$ of these labels.
In this case we see that the Pfaffian is nonzero if and only if
$$\#\{j\in J~|~j\equiv 1\bmod 2\}$$ is even.   Since $(n-a')/2$ of the labels are even and $(n-a')/2$ are odd, arguing as in the third case, we see that  
\begin{eqnarray*}
f_{0,0}(n) &=& \sum_{a'+b'\le n} \sum_{i=0}^{b'} {(n-a')/2\choose i}{(n-a')/2\choose b'-i}\chi(i \equiv b'\equiv a'-1\equiv 0\bmod 2) \\
&=& \sum_{j=0}^{n/2} \sum_{i=0}^{n/2-j}\sum_{k=0}^{n/2-j} {n/2 - j\choose i}{n/2-j+1\choose k}\chi(i\equiv k\equiv 0\bmod 2) \\
&=& 1+\sum_{j=0}^{n/2-1} 2^{n/2-j-1}2^{n/2-j-1} \\
&=& 1+(1+2^2 +\cdots + 2^{n-2}). \end{eqnarray*}
Now let $f(n)$ denote the number of Cauchon diagrams in $\mathcal{C}'_{2,n}$ with nonzero Pfaffian. 
Then we see that if $n$ is odd, $$f(n)=2^{n}-1$$ and if $n\ge 2$ is even then
$$f(n)=f_{0,0}(n)+f_{0,1}(n)=1+(1+2+4+\cdots + 2^{n-2}) = 2^{n-1}.$$
We now put this information together to obtain the desired result.  By 
Theorem \ref{thm: prim}, the number of primitive $\hc$-primes is just the number of $2\times n$ Cauchon diagrams with nonzero Pfaffian.  Since adding a completely black column does not affect the Pfaffian, we see that for $n\ge 1$ this number is just
\begin{eqnarray*}
&~& 1+\sum_{0<m\le n} {n\choose m} f(m)\\
&=& 1+\sum_{m\le n} {n\choose m}(2^m-1)\chi(m\equiv 1\bmod 2) + \sum_{0<m\le n}{n\choose m}2^{m-1}\chi(m\equiv 0\bmod 2) \\
&=& 1+\frac{1}{2}\Bigg( (2+1)^n - (1-2)^n-2^{n}\Bigg) + \frac{1}{4}\Bigg((1+2)^n+(1-2)^n-2\Bigg) \\
&=& \frac{3^{n+1}-2^{n+1}+(-1)^{n+1}+2}{4}. \end{eqnarray*}
This completes the proof. 
\end{proof}
\begin{cor} Then the proportion of  primitive $\hc$-primes in $\oqm2n$ tends to $3/8$ as $n\rightarrow \infty$.
\end{cor}
\begin{proof} We have just shown that number of primitive $\hc$-primes in $\oqm2n$ is asymptotic to $3^{n+1}/4$ as $n\rightarrow \infty$.  On the other hand, as the number of $\hc$-primes in $\oqm2n$ is equal to the number of $\hc$-primes in $\co_q(M_{n,2})$, it follows from \cite[Corollary 1.5]{laucomb} that the total number of $\hc$-primes in $\oqm2n$ is equal to $2 \cdot 3^n - 2^n$. Hence the proportion of primitive $\hc$-primes in $\oqm2n$ is asymptotic to $3/8$.  
\end{proof}
\subsection{Data and conjectures.}
Let $P(m,n)$ denote the number of primitive $\hc$-prime ideals in $\oqmmn$.  Using Maple, we obtained the following data.  
\begin{table}[ht]
\caption{The values of $P(m,n)$ for small $m$ and $n$}
\centering
\begin{tabular}{|r|r|r|r|r|r|r|r|r|r|}
\hline
$m$ & $P(m,1)$ & $P(m,2)$ & $P(m,3)$ & $P(m,4)$ & $P(m,5)$ & $P(m,6)$ & $P(m,7)$ & $P(m,8)$ & $P(m,9)$\\
\hline
1 & 1 & 2 &   4 &   8      & 16      &   32     & 64      & 128      & 256       \\
2 & 2 & 5 &   17 &   53      & 167      &   515     & 1577      & 4793      & 14507      \\
3 & 4 & 17 & 70  & 329 & 1414 & 6167 & 25960 & 108629 & 447874  \\
4 & 8 & 53 & 329 & 1865 & 11243 &&&&\\
5 & 16 &167 & 1414 & 11243 &  80806  &&&&\\
\hline 
\end{tabular}
\end{table}

We know formulas for $P(1,n)$ and $P(2,n)$ and so it is natural ask if this can be extended.  We thus pose the following question.
\begin{quest}
Can a closed formula for $P(m,n)$ be given?  In particular, can a closed formula for the diagonal terms, $P(n,n)$, be given?
\end{quest}
Using this table and the analogy with the $1\times n$ and $2\times n$ cases, we make the following conjecture.
\begin{conj} The number of primitive $\hc$-primes in the ring $\co_q(M_{3,n})$ is given by $$\frac{1}{8} \cdot \left( 15\cdot 4^n - 18 \cdot 3^n +13 \cdot 2^n - 6\cdot(-1)^n + 3\cdot (-2)^n\right).$$
\end{conj}
More generally we believe the following holds.
\begin{conj} Let $m\ge 1$ be a positive integer and let $P(m,n)$ denote the number of primitive $\hc$-prime ideals in $\oqmmn$.   Then there exist rational constants $c_{m+1},c_{m},\ldots ,c_{2-m}$ such that
$$P(m,n) \ = \ \sum_{j=2-m}^{m+1} c_j j^n$$ for all positive integers $n$.  
Moreover, $c_{m+1}=1\cdot 3 \cdot 5 \cdots (2m-1)/2^m$.
\end{conj}
This conjecture, if true, would imply the truth of the following conjecture.
\begin{conj} Let $m$ be a fixed positive integer.  Then the proportion of $\hc$-primes in $\oqmmn$ that are primitive tends to ${2m \choose m}/4^m$ as $n\rightarrow\infty$.
\end{conj} 

We note that if one follows the Proof of Theorem \ref{thm: criterion}, then one sees that the Pfaffian of a $2\times n$ Cauchon diagram is always either $0$ or plus or minus a power of $2$.  This also appears to be the case for larger Cauchon diagrams.  We therefore make the following conjecture.
\begin{conj} Let $C$ be a labeled $m\times n$ Cauchon diagram.  Then
$|{\rm Pfaffian}(C)|$ is either $0$ or a power of $2$.  
\end{conj}
We note this conjecture, if true, would allow us to simplify many computations since to determine if the 
Pfaffian is nonzero, it would suffice to consider it mod $3$.\\

\begin{flushleft}
\textbf{Acknowledgments.} We thank Ken Goodearl and the anonymous referee for useful comments on a previous draft of this paper. Also, the second author would like to thank Lionel Richard for interesting conversations on the topics of this paper. Part of this work was done while the second author was visiting Simon Fraser University. He wishes to thank NSERC for supporting his visit.
\end{flushleft}


\vskip 1cm

\noindent J. Bell:\\
Department of Mathematics, Simon Fraser University,\\
8888 University Drive,\\
Burnaby, BC, V5A 1S6, Canada\\
E-mail : jpb@math.sfu.ca
\\

\noindent S. Launois:\\
Institute of Mathematics, Statistics and Actuarial Science,\\
University of Kent at Canterbury, CT2 7NF, UK\\
E-mail : S.Launois@kent.ac.uk
\\

\noindent N. Nguyen:\\
Department of Mathematics, Simon Fraser University,\\
8888 University Drive,\\
Burnaby, BC, V5A 1S6, Canada\\
E-mail : tnn@sfu.ca
\\

 \end{document}